\documentclass[12pt,leqno]{article}

\usepackage[cp1250]{inputenc}
\usepackage[OT4]{fontenc}
\usepackage{amsfonts}
\usepackage{amsmath}
\usepackage{amssymb}
\usepackage{amsthm}
\usepackage{fancyhdr}
\newcommand{\R}{\mathbb{R}}

\newcommand{\codim}{\mathop{\rm codim}\nolimits}
\newcommand{\rank}{\mathop{\rm rank}\nolimits}
\newcommand{\corank}{\mathop{\rm corank}\nolimits}

\newcommand{\sgn}{\mathop{\rm sgn}\nolimits}
\newcommand{\ar}{\longrightarrow}
\newcommand{\inv}{^{-1}}

\newtheorem{theorem}{Theorem}[section]
\newtheorem{lemma}[theorem]{Lemma}

\newtheorem{prop}[theorem] {Proposition}

\theoremstyle{definition}
\newtheorem{definition}[theorem]{Definition}

\newtheorem{ex}[theorem]{Example}

\title{\thanks{%
Iwona Krzy\.{z}anowska and Aleksandra~Nowel\\
University of Gda\'{n}sk,
              Institute of Mathematics \\
              80-952 Gda\'{n}sk, Wita Stwosza 57\\
              Poland\\
              Tel.: +48-58-5232059\\
              Fax: +48-58-3414914\\
              Email: Iwona.Krzyzanowska@mat.ug.edu.pl\\
              Email: Aleksandra.Nowel@mat.ug.edu.pl\\ \\
{\em Keywords:} singularity, fold, cusp, one--generic, quadratic form\\              
2010 \emph{Mathematics Subject Classification} primary: 58K05, secondary: 14Q20, 14P99}
Criteria for singularities for mappings from two--manifold to the plane. The number and signs of cusps.
}

\author{Iwona~Krzy\.{z}anowska \and Aleksandra~Nowel}

\date{2016}

\begin{document}

\def\nothanksmarks{\def\thanks##1{\protect\footnotetext[0]{\kern-\bibindent##1}}}

\nothanksmarks

\maketitle

\pagestyle{fancy}

\lhead{\fancyplain{}{\textsc{\small Krzy\.{z}anowska, Nowel}}}
\rhead{\fancyplain{}{\emph{\small Criteria for singularities. Number and signs of cusps.}}}

\begin{abstract}
Let $M\subset \R^{n+2}$ be a two--dimensional complete intersection. We show how to check whether a mapping
$f\colon M\longrightarrow\R^2$ is $1$--generic with only folds and cusps as singularities. In this case we give an effective method to count the number of positive and negative cusps of a polynomial $f$, using the signatures of some quadratic forms.
\end{abstract}

\section{Introduction}

In \cite{whitney}, Whitney investigated a smooth mapping between two surfaces. He proved that for a generic mapping the only possible types of singular points are folds and simple cusps. With smooth oriented $2$--dimensional manifolds $M$ and $N$, and a smooth mapping $f\colon M\to N$ with a simple cusp $p\in M$ one can associate a sign $\mu (p)=\pm 1$ defined as the local topological degree of the germ of $f$ at $p$.

In \cite{krzysza}, the authors studied smooth mappings from the plane to the plane, and they presented methods of checking whether a map is a generic one with only folds and simple cusps as singular points. They also gave the effective formulas to determine the number of positive and negative cusps in therms of signatures of quadratic forms.

Criteria for types of Morin singularities of mappins from $\R^m$ to $\R^n$ (in case $m=n=2$ they are folds and cusps) were presented in \cite{saji}. Moreover some results concerning the algebraic sum of cusps are contained in \cite{fukudaishikawa}, \cite{quine}, and in \cite{farjelrua} in the complex case.

In this paper we investigate properties of mappings $f=\tilde{f}|_M \colon M\to \R^2$, where $M=h^{-1}(0)$ is a $2$--dimensional complete intersection, $h\colon \R^{n+2}\to \R^{n}$, $\tilde{f}\colon \R^{n+2}\to \R^2$. We give methods for checking whether $f$ is $1$--generic (in sense of \cite{golub}) and whether a given singular point $p\in M$ of $f$ is a fold point or a simple cusp (Theorem \ref{one-generic}, Propositions \ref{fold}, \ref{simple-cusp}). We define $F\colon \R^{n+2}\to \R^2$ associated with $\tilde{f}$ and $h$ such that for a simple cusp $p$ of $f$ the sign of it $\mu(p)=\sgn \det \left [
\begin{matrix}
DF(p)\\
Dh(p)
\end{matrix}
\right ]$ (Theorem \ref{signofcusp}).

In the case where $\tilde{f}$ and $h$ are polynomial mappings, we construct an ideal $S\subset \R[x]=\R[x_1,\ldots ,x_{n+2}]$ such that if $S=\R[x]$ then $f$ is $1$--generic with only folds and simple cusps as singular points (Proposition \ref{effective-sing}). Then we define an ideal $J$ such that the set of its real zeros $V(J)$ is the set of simple cusps of $f$. If $S=\R[x]$ and $\dim _{\R} \R[x]/J<\infty$ then the number of simple cusps and the algebraic sum of them can be expressed in terms of signatures of some associated quadratic forms (Proposition \ref{effective-signs}).

\section{Preliminaries}

Let $M,N$ be smooth manifolds such that $m=\dim M$ and $n=\dim N$. Take $p\in M$. For smooth mappings $f,g:M\longrightarrow N$ such that $f(p)=g(p)=q$, we say that $f$ has first order contact with $g$ at $p$ if $Df(p)=Dg(p)$, as mappings $T_pM\longrightarrow T_qN$. Then $J^1(M,N)_{(p,q)}$ denotes a set of equivalence classes of mappings $f:M\longrightarrow N$, where $f(p)=q$, having the same first order contact at $p$. Let $$J^1(M,N)=\bigcup_{(p,q)\in M\times N} J^1(M,N)_{(p,q)}$$ denote the $1$--jet bundle of smooth mappings from $M$ to $N$.

With any smooth $f:M\longrightarrow N$ we can associate a canonical mapping $j^1f:M\longrightarrow J^1(M,N)$. Take $\sigma \in J^1(M,N)$, represented by $f$. Then by $\corank \sigma$ we denote the $\corank Df(p)$. Put $S_r=\{\sigma\in J^1(M,N)\ |\ \corank \sigma =r\}$. According to \cite[II, Theorem 5.4]{golub}, $S_r$ is a submanifold of $J^1(M,N)$, with $\codim S_r=r(|m-n|+r)$. Put $S_r(f)=\{x\in M\ |\ \corank Df(p)=r\}=(j^1f)\inv (S_r)$.

\begin{definition}
We say that $f:M\longrightarrow N$ is $1$--generic if $j^1f\pitchfork S_r$, for all $r$.
\end{definition}

According to \cite[II, Theorem 4.4]{golub}, if $j^1f\pitchfork S_r$ then either $S_r(f)=\emptyset$ or $S_r(f)$ is a submanifold of $M$, with $\codim S_r(f)=\codim S_r$.

\bigskip

In the remaining we will need the following useful fact.

\begin{lemma}\label{submersions}
Let $M$, $N$ and $P$ be smooth manifolds, and let $f\colon M\to N$, $a\colon P \to M$, $b\colon P\to N$ be such that $b=f\circ a$. If $a$ is a surjective submersion, $b$ is smooth, then $f$ is also smooth. If in addition $b$ is a submersion, then so is $f$. 
\end{lemma}

\bigskip

Let
\[h=(h_1,\ldots,h_n) :\R ^{n+k} \ar \R ^{n}\]
\[f=(f_1,\ldots,f_l) :\R^ {n+k} \ar \R ^{l}\]
be $C^1$ maps, $M:=h^{-1}(0)$. Suppose that each point $p\in M$ is a regular point of $h$, i. e.
$\rank Dh(p)=n$ in each $p\in M$. Then $M$ is an orientable $C^1$ $k$--manifold called a complete intersection. It is easy to verify that for each point $p\in M$
\begin{equation} \label{rank}
\rank Df|_M(p)=\rank
\left [
\begin{array}{c}
   Df(p)\\
   Dh(p)
\end{array} \right ]
-n.
\end{equation}

\bigskip
Assume that $N=\R^2$ and $M=h^{-1}(0)$, where $h\colon \R^{n+2}\to \R^n$ is a smooth mapping such that $\rank Dh(x)=n$ for all $x\in M$. 
In that case $M$ is a smooth $2$--manifold.

We have $J^1(\R^{n+2},\R^2)\simeq \R^{n+2}\times \R^2\times M(2,{n+2})$, where $M(2,{n+2})$ is the space of real $2\times (n+2)$--matrices.

Let us define 
\[
G=\{\sigma=(x,y,A)\in J^1(\R^{n+2},\R^2)\ |\ x\in M\}=\bigcup_{(p,q)\in M\times \R^2} J^1(\R^{n+2},\R^2)_{(p,q)}. 
\]
Then $G$ is a submanifold of $J^1(\R^{n+2},\R^2)$, and $\dim G=2n+8$.

We define a relation $\sim$ in $G$:
$(x_1,y_1,A_1)=\sigma _1\sim \sigma _2=(x_2,y_2,A_2)$ if and only if $x_1=x_2$ and $y_1=y_2$, and $A_1|_{T_{x_1}M}=A_2|_{T_{x_1}M}$ considered as linear mappings on $T_{x_1}M\subset T_{x_1}\R^{n+2}$.

\begin{prop}
$G/_{\sim}$ is a smooth manifold diffeomorphic to $J^1(M,\R^2)$ such that the projection $pr\colon G\to G/_{\sim}$ is a submersion. 
\end{prop}

\begin{proof}
Using \cite[Part II, Chap. III, Sec. 12, Th. 1 and Th. 2]{serre}, to verify that $G/_{\sim}$ is a smooth manifold such that the projection $pr\colon G\to G/_{\sim}$ is a submersion, it is enough to show that
\begin{itemize}
\item[a)] the set $R=\{(\sigma _1,\sigma _2)\in G\times G\ |\ \sigma _1\sim \sigma _2\}$ is a submanifold of $G\times G$, 
\item[b)] the projection $\pi \colon R \to G$ is a submersion. 
\end{itemize}
Take $x\in M$, then in a neighbourhood of $x$ in $\R^{n+2}$ there exists a smooth non--vanishing vector field $(v_1,v_2)\in \R^{n+2}\times \R^{n+2}$ such that 
\[
\operatorname{Span} \{v_1,v_2\}=\left ( \operatorname{Span} \{\nabla h_1,\ldots ,\nabla h_n\}\right ) ^{\bot} 
\]
at every point of this neighbourhood. Then at points of $M$ vectors $v_1,v_2$ span the tangent space to $M$.  

Let us define $\gamma \colon J^1(\R^{n+2},\R^2)\times J^1(\R^{n+2},\R^2)\to \R^{2n+8}$ by
\[
\gamma (\sigma_1,\sigma_2)=\gamma((x_1,y_1,A_1),(x_2,y_2,A_2))=
\]
\[
=(x_1-x_2,y_1-y_2,A_1v_1(x_1)-A_2v_1(x_1),A_1v_2(x_1)-A_2v_2(x_1),h(x_1)). 
\]
Hence $\gamma (\sigma_1,\sigma_2)=0$ if and only if $(\sigma_1,\sigma_2)\in R$. Then locally $\gamma ^{-1}(0)=R$. Moreover $\gamma$ is a submersion at points from $R$, so $R$ is a submanifold of $G\times G$, and a) is proven.

Using equation (\ref{rank}) it is easy to see that $\rank D\pi=2n+8=\dim G$, so $\pi$ is a submersion and we have b).

Now we will prove that $G/_{\sim}$ is diffeomorphic to $J^1(M,\R^2)$. Since $M$ is a submanifold of $\R^{n+2}$, there exists a tubular neighbourhood $U$ of $M$ in $\R^{n+2}$ with a smooth retraction $r\colon U\to M$, which is also a submersion.

Let us define $\Psi\colon J^1(M,\R^2)\to G/_{\sim}$ by
\[
\Psi(\sigma)=\Psi ([g])=[g\circ r]\in G/_{\sim}. 
\]
Note that $\Psi$ is a well--defined bijection and $\Psi ^{-1}$ is given by $G/_{\sim}\ni [g]\mapsto [g|_{M}]\in J^1(M,\R^2)$. The mapping $\Psi ^{-1}\circ pr\colon G\to J^1(M,\R^2)$ can be given by $G\ni [g]\mapsto [g|_{M}]\in J^1(M,\R^2)$ and we see that it is a smooth submersion. So according to Lemma \ref{submersions}, $\Psi^{-1}$ is also a smooth submersion. Since $\Psi^{-1}$ is bijective, it is a diffeomorphism.
\end{proof}

\section{Checking 1--genericity and recognizing folds and cusps}

Let $\tilde{f}\colon \R^{n+2}\to \R^2$ be smooth and put $f=\tilde{f}|_M\colon M\to \R^2$, where $M=h^{-1}(0)$ is a $2$--dimensional complete intersection. Using mappings $h$ and $\tilde{f}$ defined on $\R^{n+2}$,
we will present an effective method to check whether $f$ is $1$--generic.

Put $\Phi\colon G/_{\sim}\to \R$ as
\[
\Phi([(x,y,A)])=\det 
\left [
\begin{matrix}
A \\
Dh(x)
\end{matrix}
\right ].
\]
Notice that if $[(x,y,A)]\in G/_{\sim}$ is represented by $g$ defined near $x\in \R^{n+2}$, then $\Phi([g])=\det 
\left [
\begin{matrix}
Dg(x) \\
Dh(x)
\end{matrix}
\right ]$.  

\begin{lemma}
$\Phi$ is well--defined. 
\end{lemma}

\begin{proof}
Take $(x,y,A_1)$ and $(x,y,A_2)$ representing the same element in $G/_{\sim}$. Then $A_1v_1=A_2v_1$ and $A_1v_2=A_2v_2$, where $v_1,v_2\in \R^{n+2}$ span $T_xM$, and so they both are orthogonal to all vectors $\nabla h_i(x)$. 

Hence we have
\[
\det \left (
\left [
\begin{matrix}
A_1\\
Dh(x)
\end{matrix}
\right ]
\left [
\begin{matrix}
v_1 & v_2 & \nabla h_1(x) & \ldots & \nabla h_n(x)
\end{matrix}
\right ]\right )=
\det
\left [
\begin{matrix}
A_1v_1\ A_1v_2 & *\\
\mathbf{0} & Dh(x)Dh(x)^T
\end{matrix}
\right ]=
\]
\[
\det
\left [
\begin{matrix}
A_2v_1\ A_2v_2 & **\\
\mathbf{0} & Dh(x)Dh(x)^T
\end{matrix}
\right ]=
\det \left (
\left [
\begin{matrix}
A_2\\
Dh(x)
\end{matrix}
\right ]
\left [
\begin{matrix}
v_1 & v_2 & \nabla h_1(x) & \ldots & \nabla h_n(x)
\end{matrix}
\right ]\right ).
\]
Since $\det [v_1 \ v_2 \ \nabla h_1(x) \ \ldots \ \nabla h_n(x)]\neq 0$, we obtain
\[
\det
\left [
\begin{matrix}
A_1\\
Dh(x)
\end{matrix}
\right ]=
\det
\left [
\begin{matrix}
A_2\\
Dh(x)
\end{matrix}
\right ].
\]
\end{proof}

\begin{lemma} \label{rank12}
$\Phi$ is a submersion at every $[(x,y,A)]\in G/_{\sim}$ such that $\rank  \left [
\begin{matrix}
A\\
Dh(x)
\end{matrix}
\right ]\geqslant n+1$. 
\end{lemma}

\begin{proof}
Put $\tilde{\Phi}\colon G \to \R$ as $\tilde{\Phi}(x,y,A)=\det \left [
\begin{matrix}
A\\
Dh(x)
\end{matrix}
\right ]$. Then $\tilde{\Phi}(x,y,A)$ can be expressed as a linear combination of elements of one of rows of the matrix $A$, whose coefficients are appropriates $(n+1)$--minors of the matrix $\left [
\begin{matrix}
A\\
Dh(x)
\end{matrix}
\right ]$. Since at least one of these minors is not $0$, $\tilde{\Phi}$ is a submersion at $(x,y,A)$. Notice that $\tilde{\Phi}=\Phi \circ pr$, so by Lemma \ref{submersions}, $\Phi$ is a submersion at $[(x,y,A)]$.  
\end{proof}

For a smooth mapping $\tilde{f}\colon \R^{n+2}\to \R^2$ we define $d\colon \R^{n+2}\to \R$ as
\[
d(x)= \det \left [
\begin{matrix}
D\tilde{f}(x)\\
Dh(x)
\end{matrix}
\right ].
\]

According to (\ref{rank}) for $f=\tilde{f}|_M\colon M\to \R^2$ we have $x\in S_i(f)$ if and only if $\rank \left [
\begin{matrix}
D\tilde{f}(x)\\
Dh(x)
\end{matrix}
\right ]=n+2-i$, for $i=1,2$, and so $S_1(f)\cup S_2(f)=d^{-1}(0)\cap M$.

\begin{theorem} \label{one-generic}
A mapping $f=\tilde{f}|_{M}\colon M \to \R^2$ is $1$--generic if and only if 
$d|_M$ is a submersion at points from 
$d^{-1}(0)\cap M$,
i. e. $\rank \left [
\begin{matrix}
Dd(x)\\
Dh(x)
\end{matrix}
\right ]=n+1$, 
for $x\in d^{-1}(0)\cap M$.
If that is the case, then $S_1(f)=d^{-1}(0)\cap M$.
\end{theorem}

\begin{proof}
Let $x\in S_1(f)$. According to Lemma \ref{rank12}, $\Phi$ is a submersion at $\Psi (j^1f(x))$. Notice that there exists a small enough neighbourhood $U$ of $\Psi (j^1f(x))$ such that $\Phi|_U$ is a submersion and 
\[
U\cap \Psi(S_1)={\Phi|_U}^{-1}(0). 
\]
We have $j^1f\pitchfork S^1$ at $x$ if and only if $\Psi(j^1f)\pitchfork \Psi(S^1)$ at $x$.
According to \cite[II, Lemma 4.3]{golub}, $\Psi(j^1f)\pitchfork \Psi(S^1)$ at $x$ if and only if $\Phi|_U\circ \Psi\circ j^1f$ is a submersion at $x$.

Let us see that $\Phi|_U\circ \Psi\circ j^1f(x)=d(x)$ for $x\in M$. We get that for $x\in S_1(f)$, $j^1f\pitchfork S^1$ at $x$ if and only if $d|_M\colon M\to \R$ is a submersion at $x$, i. e. $\rank \left [
\begin{matrix}
Dd(x)\\
Dh(x)
\end{matrix}
\right ]=n+1$.

Note that since $\codim S_2=4$, $j^1f \pitchfork S_2$ if and only if $S_2(f)=\emptyset$. On the other hand, if $x\in S_2(f)$, then 
\[
\rank \left [
\begin{matrix}
D\tilde{f}(x)\\
Dh(x)
\end{matrix}
\right ]=n,
\]
the elements of $Dd(x)=D \left (\det \left [
\begin{matrix}
D\tilde{f}(x)\\
Dh(x)
\end{matrix}
\right ]  \right)$ are linear combinations of $(n+1)$--minors of this matrix,
and so $Dd(x)=(0,\ldots ,0)$. We get that if $d|_M$ is a submersion at points from 
$d^{-1}(0)\cap M$, then $S_2(f)=\emptyset$.
\end{proof}

From now on we assume that $f=\tilde{f}|_M\colon M\to \R^2$ is $1$--generic. Then by Theorem \ref{one-generic}, for $x$ near $S_1(f)$, the vectors $\nabla h_1(x),\ldots ,\nabla h_n(x),\nabla d(x)$ are linearly independent and $S_1(f)$ is $1$--dimensional submanifold of $M$.

For $x\in \R^{n+2}$ and the matrix $\left [
\begin{matrix}
Dd(x)\\
Dh(x)
\end{matrix}
\right ]$, by $w_i(x)$ we will denote its $(n+1)$--minors obtained by removing $i$--th column. We define a vector field $v\colon \R^{n+2}\to \R^{n+2}$ as
\[
v(x)=\left(-w_1(x),w_2(x),\ldots ,(-1)^{n+2}w_{n+2}(x)\right ). 
\]
Then for $x\in S_1(f)$ the vector $v(x)$ is a generator of 
\[
T_xS_1(f)=\left (\operatorname{Span} \{\nabla h_1(x),\ldots ,\nabla h_n(x),\nabla d(x)\}\right )^{\bot}. 
\]
Put $F=(F_1,F_2)\colon \R^{n+2}\to \R^2$ as
\[
F(x)=D\tilde{f}(x)(v(x)). 
\]

We will call $p\in S_1(f)$ a \textbf{fold point} if it is a regular point of $f|_{S_1(f)}$. 

\begin{prop}\label{fold}
For a $1$--generic $f$ and a point $p\in S_1(f)$ the following are equivalent:
\begin{itemize}
\item[(a)] $p$ is a fold point; 
\item[(b)] $\rank \left [
\begin{matrix}
D\tilde{f}(p)\\
Dh(p)\\
Dd(p)
\end{matrix}
\right ]=n+2$; 
\item[(c)] $F(p)\neq 0$.
\end{itemize}
\end{prop}
\begin{proof}
Since $f$ is $1$--generic, $S_1(f)=(h,d)^{-1}(0)$ is a complete intersection, and so
the equivalence of the first two conditions is a simple consequence of the equation (\ref{rank}).   

We see that $F(p)\neq 0$ iff $\langle \nabla \tilde{f}_1(p), v(p)\rangle\neq 0$ or $\langle \nabla \tilde{f}_2(p), v(p)\rangle\neq 0$ iff at least one of $\nabla \tilde{f}_1(p), \nabla \tilde{f}_2(p)$ does not belong to $\operatorname{Span} \{\nabla h_1(x),\ldots ,\nabla h_n(x),\nabla d(x)\}$ iff $\rank \left [
\begin{matrix}
D\tilde{f}(p)\\
Dh(p)\\
Dd(p)
\end{matrix}
\right ]=n+2$. So we get (b) $\Leftrightarrow$ (c).
\end{proof}

If $f=(f_1,f_2)\colon M\to \R^2$ is $1$--generic, then for $p\in S_1(f)$ one of the following two conditions can occur. \begin{equation}\label{condition1}
T_pS_1(f)+ \ker Df(p)=\R^2 ,
\end{equation} 
\begin{equation}\label{condition2}
T_pS_1(f)=\ker Df(p).
\end{equation}
It is easy to see that $p\in S_1(f)$ satisfies  (\ref{condition1}) if and only if $F(p)\neq 0$, and then $p$ is a fold point. 

Assume that  condition (\ref{condition2}) holds at $p\in S_1(f)$. By the previous Proposition this is equivalent to the condition $F(p)=0$.

Take  a smooth function $k$  on $M$ such that $k\equiv 0$ on $S_1(f)$ and $Dk(p)\neq 0$
(our mapping $d|_M$ satisfies both these conditions).
Let $\xi$ be a non--vanishing vector field  along $S_1(f)$ such that  $\xi $ is in the kernel of $Df$
at each point of $S_1(f)$ near $p$. Then $Dk(\xi )$ is a function on $S_1(f)$ having a zero at $p$.
The order of this zero does not depend on the choice of $\xi$ or $k$ (see \cite[p. 146]{golub}),
so in our case it equals the order of $Dd|_M(\xi)$ at $p$. Following \cite{golub} we will say that $p$ is a \textbf{simple cusp} (or \textbf{cusp} for short) if $p$ is a simple zero of $Dd|_M(\xi)$. If this is the case, then locally near $p$ the mapping $f$ has a form $(x_1,x_2)\mapsto (x_1,x_2^3+x_1x_2)$ (see \cite{whitney}, \cite{golub}).

\begin{prop} \label{simple-cusp}
Assume that $f$ is $1$--generic and $p\in S_1(f)$. Then $p$ is a simple cusp if and only if $F(p)=0$ and $\rank \left [
\begin{matrix}
DF(p)\\
Dh(p)\\
Dd(p)
\end{matrix}
\right ]=n+2$. 
\end{prop}
\begin{proof}
Take $p\in S_1(f)$. Note that $F(p)=0$ is equivalent to the condition $T_pS_1(f)=\ker Df(p)$. So we assume that $F(p)=0$.

Let us take a small neighbourhood $U\subset \R^{n+2}$ of $p$ and a smooth vector field $w\colon U\to \R^{n+2}$ such that 
\[
\operatorname{Span} \{ w(x)\} = \left (\operatorname{Span} \{\nabla h_1(x),\ldots ,\nabla h_n(x), v(x)\}\right )^{\bot} 
\mbox{ and }
\langle \nabla d(x),w(x)\rangle \neq 0, 
\]
for $x\in U$.
We define a smooth vector field $\xi _i\colon S_1(f)\cap U\to \R^{n+2}$ for $i=1,2$ by
\[
\xi_i(x)=\frac{F_i(x)}{\langle \nabla d(x),w(x)\rangle}w(x)-\frac{\langle \nabla \tilde{f}_i(x),w(x)\rangle}{\langle \nabla d(x),w(x)\rangle}v(x). 
\]
By our assumptions
\[
\rank \left [
\begin{matrix}
D\tilde{f}(p)\\
Dh(p)
\end{matrix}
\right ]= 
\rank \left [
\begin{matrix}
Dd(p)\\
Dh(p)
\end{matrix}
\right ]=\rank \left [
\begin{matrix}
D\tilde{f}(p)\\
Dd(p)\\
Dh(p)
\end{matrix}
\right ]=n+1,
\]
and then there exist $\alpha, \beta\in \R$ such that $\alpha ^2+\beta ^2\neq 0$,
$\nabla d(p)=\alpha \nabla \tilde{f}_1(p)+\beta \nabla \tilde{f}_2(p)+$ some linear combination of $\nabla h_i(p)$.
So 
\[
0\neq \langle \nabla d(p),w(p)\rangle=\alpha \langle \nabla \tilde{f}_1(p),w(p)\rangle+\beta  \langle \nabla \tilde{f}_2(p),w(p)\rangle,
\]
and then $\langle \nabla \tilde{f}_1(p),w(p)\rangle\neq 0$ or $\langle \nabla \tilde{f}_2(p),w(p)\rangle\neq 0$. Hence at least one of $\xi_i(p)=-\frac{\langle \nabla \tilde{f}_i(p),w(p)\rangle}{\langle \nabla d(p),w(p)\rangle}v(p)$ is different from $0$.
Of course $\xi_i(p)\in T_pS_1(f)=\operatorname{Span}\{v(p)\}$.

Since for $x\in S_1(f)\cap U$ we have $\xi_i(x)\in \left (\operatorname{Span} \{\nabla h_1(x),\ldots ,\nabla h_n(x)\}\right )^{\bot}$, $\langle \nabla \tilde{f}_i(x),\xi_i(x)\rangle=0$, and $\rank \left [
\begin{matrix}
D\tilde{f}(x)\\
Dh(x)
\end{matrix}
\right ]=n+1$.
It is easy to see that 
\[
\left [
\begin{matrix}
D\tilde{f}(x)\\
Dh(x)
\end{matrix}
\right ]\xi_i(x)=0, 
\]
and so $\xi_i(x)\in \ker(Df(x))$ for $i=1,2$.

Notice that $Dd|_M(x)\xi_i(x)=\langle \nabla d(x),\xi_i(x)\rangle=F_i(x)$ for $x\in S_1(f)\cap U$. Take $i$ such that $\xi_i(p)\neq 0$. We get that $p$ is a simple cusp if and only if $p$ is a simple zero of $F_i|_{S_1(f)}$, then $\rank \left [
\begin{matrix}
DF(p)\\
Dh(p)\\
Dd(p)
\end{matrix}
\right ]=n+2$.

On the other hand, if for $j=1,2$, $\rank \left [
\begin{matrix}
DF_j(p)\\
Dh(p)\\
Dd(p)
\end{matrix}
\right ]=n+2$, then $p$ is a simple zero of $F_i|_{S_1(f)}$. So let us assume, that for example $\rank \left [
\begin{matrix}
DF_2(p)\\
Dh(p)\\
Dd(p)
\end{matrix}
\right ]=n+1$
and $\rank \left [
\begin{matrix}
DF_1(p)\\
Dh(p)\\
Dd(p)
\end{matrix}
\right ]=n+2$. Since for $x\in S_1(f)\cap U$, $\rank \left [
\begin{matrix}
D\tilde{f}(x)\\
Dh(x)
\end{matrix}
\right ]=n+1$, there exist smooth $\alpha, \beta$ such that
$\alpha^2(x)+\beta^2(x)\neq 0$ and $\alpha(x)F_1(x)+\beta(x)F_2(x)=0$
for $x\in S_1(f)\cap U$. Then differentiating the above equality in $S_1(f)\cap U$ we get $\beta (p)\neq 0$ and we obtain $\langle \nabla\tilde{f}_2(p),w(p)\rangle=0$. So $\xi_2(p)=0$, that means $i$ must be $1$, and $\rank \left [
\begin{matrix}
DF_i(p)\\
Dh(p)\\
Dd(p)
\end{matrix}
\right ]=n+2$ implies that $p$ is a simple zero of $F_i|_{S_1(f)}$.  
\end{proof}

\section{Signs of cusps}

Let $f\colon M\to \R^2$ be a smooth map on a smooth oriented $2$--dimensional manifold. For a simple cusp $p$ of $f$ we denote by $\mu(p)$ the local topological degree $\deg _pf$ of the germ $f\colon (M,p)\to (\R^2,f(p))$. From the local form of $f$ near $p$ it is easy to see that $\mu (p)=\pm 1$. We will call it the \textbf{sign of the cusp} $p$. 

In \cite{krzysza}, the authors investigated the algebraic sum of cusps of a $1$--generic mapping $g=(g_1,g_2)\colon \R^2\to \R^2$. They defined $G\colon \R^2\to \R^2$ as $G(x)=Dg(x)\zeta(x)$, where $\zeta(x)=(\zeta_1(x),\zeta_2(x))=\left (-\frac{\partial}{\partial x_2}\det Dg(x),\frac{\partial}{\partial x_1}\det Dg(x)\right )$ is tangent to $S_1(g)$ for $x\in S_1(g)$.

According to \cite[Proposition 1]{krzysza}, for a simple cusp $q\in \R^2$ of $g$, we have $\det DG(q)\neq 0$ and $\mu (q)=\sgn \det DG(q)$.

Using the facts and proofs from \cite[Section 3.]{krzysza} it is easy to show the following.
\begin{lemma} \label{signs}
Let $\eta=(\eta_1,\eta_2)$ be a non--zero vector field on $\R^2$. Assume that in some neighbourhood of the simple cusp $q$ of $g$ there exists a smooth non--vanishing function $s$ such that on $S_1(g)$ we have $s(x)\eta(x)=\zeta(x)$. Then for $\tilde{G}(x)=Dg(x)\eta(x)$
\[
\sgn \det DG(q)=\sgn \det D\tilde{G}(q). 
\]
\end{lemma}
\begin{proof}
Following \cite[Section 3.]{krzysza} we can assume that $q=0$ and there exist $\alpha, \beta \neq 0$ such that 
\[
Dg(0)=\left [
\begin{matrix}
0 & \alpha\\
0 & 0
\end{matrix}
\right ], \quad 
\zeta(0)=\left (\beta, 0\right ), \quad \frac{\partial ^2g_2}{\partial x_1^2}(0)=0.
\]
We can take a smooth $\varphi\colon (\R,0)\to (\R,0)$ such that locally $S_1(g)=\{(t,\varphi(t)\}$. Then $\varphi '(0)=0$ and
\[
\frac{d}{dt}s(t,\varphi(t))\eta_2(t,\varphi(t))=\frac{d}{dt}\zeta_2(t,\varphi(t)), 
\]
hence $s(0)\dfrac{\partial \eta_2}{\partial x_1}(0)=\dfrac{\partial \zeta_2}{\partial x_1}(0)$. Easy computations show that $\det DG(0)=s^2(0)\det D\tilde{G}(0)$.
\end{proof}

\bigskip

Let us recall that $\tilde{f}\colon \R^{n+2}\to \R^2$ is smooth and $f=\tilde{f}|_M\colon M\to \R^2$  is $1$--generic, $M=h^{-1}(0)$ is a complete intersection. In the previous section we have defined a vector field $v\colon \R^{n+2}\to \R^{n+2}$ such that for $x\in S_1(f)$ the vector $v(x)$ spans $T_xS_1(f)$, and the mapping $F(x)=D\tilde{f}(x)v(x)$.

\begin{theorem} \label{signofcusp}
Let us assume that $p$ is a simple cusp of a $1$--generic map $f\colon M\to \R^2$, where $f=\tilde{f}|_M$ and $M=h^{-1}(0)$ is a complete intersection. Then $\mu(p)=\sgn \det \left [
\begin{matrix}
DF(p)\\
Dh(p)
\end{matrix}
\right ]$. 
\end{theorem}
\begin{proof}
We can choose a chart $\phi$ of $\R^{n+2}$ defined in some neighbourhood of $p$ such that both $\phi$ and the corresponding chart $\phi_M$ of $M$, i.e. $\phi|_M=(\phi_M,0)\colon M\to \R^2\times \{ 0\}$, preserve the orientations. Put $q=\phi_M(p)$ and take $G$ as above for the mapping $g=f\circ\phi_M^{-1}\colon (\R^2,q)\to \R^2$.

For $x\in M$  we define
$\eta=(\eta_1,\eta_2)$ as $D\phi(x)v(x)=(\eta_1(x),\eta_2(x),0,\ldots,0)$. Let $y\in \R^2$ be such that $\phi(x)=(y,0,\ldots,0)$, i. e. $\phi_M(x)=y$. 
Since $\eta(x)=\eta(\phi_M^{-1}(y))$ is a non--zero vector in the tangent space at $y$ of $\phi_M(S_1(f))=S_1(g)\subset \R^2$, as well as $\zeta (y)$, there exists a smooth non--vanishing mapping $s\colon (\R^2,q)\to \R$ such that $\zeta(y)=s(y)\eta(\phi_M^{-1}(y))$ for $y\in \phi_M(S_1(f))$.

According to \cite[Proposition 1.]{krzysza},  
\[
\mu(p)=\deg _pf=\deg _qg=\sgn \det DG(q)\neq 0. 
\]
Define $\tilde{G}(y)=Dg(y)\eta(\phi_M^{-1}(y))$. Then from Lemma \ref{signs}
\[
\sgn \det DG(q)=\sgn \det D\tilde{G}(q). 
\]
Notice that
\[
F(\phi_M^{-1}(y))=D\tilde{f}(\phi^{-1}(y,0))D\phi^{-1}(y,0)D\phi(\phi^{-1}(y,0))v(\phi^{-1}(y,0))=
\]
\[
=D(\tilde{f}\circ \phi^{-1})(y,0)(\eta(\phi^{-1}(y,0)),0)=Dg(y)(\eta(\phi_M^{-1}(y))=\tilde{G}(y). 
\]
According to \cite[Lemma 3.1.]{szafran}
\[
\sgn \det D\tilde{G}(q)=\sgn \det D(F\circ \phi_M^{-1})(q)=\sgn \det \left [
\begin{matrix}
DF(p)\\
Dh(p)
\end{matrix}
\right ]. 
\]
\end{proof}

\section{Algebraic sum of cusps of a polynomial mapping}

Now we recall a well--known fact. Take an ideal $J\subset \R[x]=\R[x_1,\ldots, x_m]$ such that the $\R$--algebra ${\cal A}=\R[x]/J$ is finitely generated over $\R$, i. e. $\dim_{\R}{\cal A}<\infty$. Denote by $V(J)$ the set of real zeros of the ideal $J$.

For $h\in{\cal A}$, we denote by $T(h)$ the trace of the $\R$--linear endomorphism
${\cal A}\ni a\mapsto h\cdot a\in{\cal A}$. Then $T:{\cal A}\rightarrow\R$ is a linear functional.
Take $\delta\in\R[x]$. Let $\Theta:{\cal A}\rightarrow\R$ be the quadratic form given by
$\Theta(a)=T(\delta\cdot a^2)$.

According to \cite{becker}, \cite{pedersenetal}, the signature $\sigma(\Theta)$
of $\Theta$ equals
\begin{equation} 
\label{sumofsigns} 
\sigma(\Theta)=\sum _{p\in V(J)} \sgn\delta(p),
\end{equation} 
and if  $\Theta$ is non-degenerate then $\delta(p)\neq 0$ for each $p\in V(J)$. 

Take polynomial mappings $\tilde{f}\colon \R^{n+2}\to \R^2$ and $h=(h_1,\ldots ,h_n)\colon \R^{n+2}\to \R^n$ such that $M=h^{-1}(0)$ is a complete intersection. 
Put $f=\tilde{f}|_M\colon M\to \R^2$. Let us recall that $d(x)=\det \left [
\begin{matrix}
D\tilde{f}(x)\\
Dh(x)
\end{matrix}
\right ]$, $v(x)=(-w_1(x),w_2(x),\ldots ,(-1)^{n+2}w_{n+2}(x))$, where $w_i(x)$ are $(n+1)$--minors obtained by removing $i$--th column from the matrix $\left [
\begin{matrix}
Dd(x)\\
Dh(x)
\end{matrix}
\right ]$, and $F(x)=D\tilde{f}(x)v(x)$.

Let us define ideals $I,S\subset \R[x]=\R[x_1,\ldots ,x_{n+2}]$ as 
\[
I=\left \langle h_1,\ldots ,h_n,d,w_1,\ldots ,w_{n+2} \right \rangle, 
\]
\[
S=\left \langle h_1,\ldots ,h_n,d,F_1,F_2,\ \det \left [
\begin{matrix}
DF_1\\
Dd\\
Dh
\end{matrix}
\right ],
\ \det \left [
\begin{matrix}
DF_2\\
Dd\\
Dh
\end{matrix}
\right ]\right \rangle. 
\]
One may check that $S\subset I$.

By Theorem \ref{one-generic} and Proposition \ref{simple-cusp} we get
\begin{prop} \label{effective-sing}
\begin{itemize}
\item[(a)] If $I=\R[x]$ then $f$ is $1$--generic. 
\item[(b)] If $S=\R[x]$ then $f$ is $1$--generic, and has only folds and simple cusps as singular points. If that is the case, then the set of simple cusps $\{x\in \R^{n+2}\ | \ h_1(x)=\ldots =h_n(x)=d(x)=F_1(x)=F_2(x)=0\}$ is an algebraic set of isolated points, so it is finite. 
\end{itemize}
\end{prop}

Let us assume that $S=\R[x]$. Put $J=\langle h_1,\ldots ,h_n,d,F_1,F_2 \rangle$, and ${\cal A}=\R[x]/J$, and assume that $\dim_{\R} {\cal A}<\infty$. Then the set of simple cusps $V(J)$ of $f$ is finite and we can count the algebraic sum of cusps, i. e. $\displaystyle \sum_{p\in V(J)} \mu(p)$. Let us define quadratic forms $\Theta_1,\Theta_2\colon {\cal A}\to \R$ by
$\Theta_1(a)=T(1\cdot a^2)$,
$\Theta_2(a)=T(\delta \cdot a^2)$, where $\delta(x)=\det \left [
\begin{matrix}
DF(x)\\
Dh(x)
\end{matrix}
\right ]$. According to the formula (\ref{sumofsigns}) and Theorem \ref{signofcusp} we get the following.

\begin{prop} \label{effective-signs} Assume that $S=\R[x]$ and $\dim_{\R} {\cal A}<\infty$. Then
\begin{itemize}
\item[(a)] $\# V(J)=\sigma (\Theta _1)$.  
\item[(b)] $\displaystyle \sum_{p\in V(J)} \mu(p)=\sigma (\Theta _2)$. 
\end{itemize} 
\end{prop}

Using previous propositions and {\sc Singular} (\cite{greueletal}) we computed the following examples.

\begin{ex}
Put $\tilde{f}=(xz^2-z^2-2z,2x^3z-y^3+z^3+3yz-z^2-y)\colon \R^3\to \R^2$ and $h=x^2+y^2+z^2-1\colon \R^3\to \R$. Then $h^{-1}(0)$ is a $2$--dimensional sphere, and the mapping $f=\tilde{f}|_{h^{-1}(0)}$ is $1$--generic, has $6$ simple cusps, $3$ of them are negative. 
\end{ex}

\begin{ex}
Put $\tilde{f}=(2xz^2-y^2+2xz,-z^3+2xy-y^2-x)\colon \R^3\to \R^2$ and $h=x^2+y^2+z^2-1\colon \R^3\to \R$. Then the mapping $f=\tilde{f}|_{h^{-1}(0)}$ is $1$--generic, has $8$ simple cusps, $6$ of them are negative. 
\end{ex}

\begin{ex}
Put $\tilde{f}=(zw-2w^2-2x,3x^3-2yz^2-yw+2zw-x)\colon \R^4\to \R^2$ and $h=(x^2+y^2-1,z^2+w^2-1)\colon \R^4\to \R^2$. Then $h^{-1}(0)$ is a $2$--dimensional torus, and the mapping $f=\tilde{f}|_{h^{-1}(0)}$ is $1$--generic, has $16$ simple cusps, $8$ of them are negative. 
\end{ex}

\begin{ex} 
Put $\tilde{f}=(3z^3+x^2-xy,2y^2z-2z^3+xy-2y^2-x)\colon \R^3\to \R^2$ and $h=x^2+y^2-z\colon \R^3\to \R$. Then $h^{-1}(0)$ is a $2$--dimensional paraboloid, and the mapping $f=\tilde{f}|_{h^{-1}(0)}$ is $1$--generic, has $3$ simple cusps, all of them are negative. 
\end{ex}

\end{document}